\documentclass{amsart}

\usepackage{amssymb}

\theoremstyle{plain}
\newtheorem{thm}{Theorem}
\newtheorem{prop}[thm]{Proposition}

\newtheorem{cor}[thm]{Corollary}

\theoremstyle{definition}
\newtheorem{defn}[thm]{Definition}

\newcommand{\eps}{\varepsilon}

\newcommand{\C}{\mathbb{C}}
\newcommand{\R}{\mathbb{R}}
\newcommand{\N}{\mathbb{N}}

\newcommand{\act}{\diamond}

\DeclareMathOperator{\id}{id}
\begin{document}

\title{Exactness from proper actions}

\author{J. Brodzki}
\email{j.brodzki@soton.ac.uk}
\author{G.A. Niblo}
\email{g.a.niblo@soton.ac.uk}
\author{N.J. Wright}
\email{n.j.wright@soton.ac.uk}

\address{School of Mathematics, University of Southampton, Highfield, SO17 1BJ}

\begin{abstract}
In this paper we show that if a discrete group $G$ acts properly isometrically on a discrete space $X$ for which the uniform Roe algebra $C_u^*(X)$ is exact then $G$ is an exact group. As a corollary, we note that if the action is cocompact then the following are equivalent: The space $X$ has Yu's property A;  $C^*_u(X)$ is exact; $C_u^*(X)$ is nuclear. 
\end{abstract}

\maketitle

\section*{Introduction}

It is a remarkable result of Yu that any discrete metric space with bounded geometry that satisfies
a F\o lner-type condition, which he called  property A, also satisfies the coarse Baum-Connes conjecture \cite{Yu}. 
The case of a countable discrete group, regarded as a coarse metric space of bounded geometry, was studied by Higson and Roe \cite{HR}, Guentner and  Kaminker \cite{GK}, and Ozawa \cite{Oz}. They showed that Yu's property A is equivalent to exactness of its reduced $C^*$-algebra  $C^*_r(G)$ and to 
nuclearity of its uniform Roe algebra $C^*_u(G)$. Furthermore, Roe showed for a discrete 
bounded geometry metric space $X$ that if $X$ has property A then $C^*_u(X)$ is nuclear (\cite{R}, Proposition 11.41). 
It is tempting to conjecture therefore that for a discrete metric space of bounded geometry, 
nuclearity of $C^*_u(X)$ is in fact equivalent to property A. This conjecture holds if $X$ admits a free cocompact action by a countable group $G$: since the 
action is free we may identify $C^*_u(G)$ with a subalgebra of $C^*_u(X)$. If we assume that
$C^*_u(X)$ is nuclear, then $C^*_u(G)$ will be exact and $G$ will have property A. The fact 
that $G$ acts cocompactly on $X$ implies that $X$ and $G$ are coarsely equivalent so $X$ 
also has property A.

The conjecture is more delicate than it 
may appear since while property A is a coarse invariant, the uniform Roe algebra is not. 
For example, for a finite space $X$ with $n$ points, $C^*_u(X)$ is the algebra of $n\times n$ matrices; 
however all finite spaces are coarsely equivalent to a point.

In this paper we address the case of the conjecture when the space admits a proper, cocompact isometric action by a countable group $G$. Generalising a theorem in coarse geometry from the case of a free action to a proper action usually requires  only a minor adjustment of the argument. This is not the case in our context. The fact that the uniform Roe algebra is not a coarse invariant is manifested in the observation that when the action is not free we no longer have the required embedding of $C^*_u(G)$ into $C^*_u(X)$.  The key idea of this paper is to replace the  proper action  on $X$ by a partially defined free action on an orbit of the original action. In outline we proceed as follows.

Given an orbit $Y$ of the action of $G$ on $X$ and a splitting $\phi$ of the natural surjection $G\rightarrow Y$ we identify the algebra $C^*_u(Y)$ with the subalgebra $C^*_u(\phi(Y))$ of $C^*_u(G)$.  We construct a free partial action of $G$ on the orbit $Y$ defined in terms of the right action of $G$ on itself and the identification of $\phi(Y)$ with $Y$.  Finally we may appeal to coarse equivalences to establish the following.


\begin{thm}
\label{mainthm}
Let $X$ be a uniformly discrete metric space, and let $G$ be a countable group acting properly by isometries on $X$. If $C^*_u(X)$ is an exact $C^*$-algebra, then $G$ is an exact group.
\end{thm}

From this we deduce the following equivalence.

\begin{cor}
If $X$ is a uniformly discrete metric space admitting a proper, cocompact isometric action by a countable group $G$ then the following are equivalent.
\begin{enumerate}
\item $X$ has property A;
\item $C^*_u(X)$ is nuclear;
\item $C^*_u(X)$ is exact;
\item $G$ is exact.
\end{enumerate}
\end{cor}

\begin{proof} {\it \ of Corollary 2}

The implication $1\Rightarrow 2$ was established by Roe  in \cite[Prop. 11.41]{R}. 
That  $2$ implies $3$ is a well known result, see, e.g.  \cite{Wassermann}. Theorem \ref{mainthm}
provides the implication $3 \Rightarrow 4$. 
Finally,  as the action of $G$ on $X$ is cocompact, $X$ and $G$ are coarsely equivalent, which gives the implication $4\, \Rightarrow \, 1$. 
\end{proof}

%

\section*{Background}
Throughout the paper we will assume that $G$ is a countable group equipped with the unique (up to coarse equivalence) proper left invariant metric $d_G$. For $R\geq 0$ we will use the notation $B_R(g)$ to denote the closed $R$-ball about $g$ in $G$.

\begin{defn}
A uniformly discrete metric space $(X,d_X)$ has \emph{property $A$} if for all $R,\eps>0$ there exists a family of finite non-empty subsets $A_x$ of $X\times \N$, indexed by $x$ in $X$, such that
\begin{itemize}
\item for all $x,y$ with $d_X(x,y)<R$ we have $\frac{|A_x\Delta A_{y}|}{|A_x\cap A_{y}|}<\eps$;
 
\item there exists $S$ such that for all $x$ and $(y,n) \in A_x$ we have $d_X(x,y)\leq S$.
\end{itemize}
\end{defn}

\begin{defn}
A kernel $u\colon X\times X \to \R$ has \emph{$(R,\eps)$-variation} if $d_X(x,y)\leq R$ implies $|u(x,y)-1|<\eps$.
\end{defn}

\begin{thm}[\cite{Tu}, Proposition 3.2]
\label{kernels}
If $X$ is a bounded geometry discrete metric space, then $X$ has property A if and only if for all $R,\eps>0$ there exists a positive definite kernel $u$ with $R,\eps$ variation, and such that there exists $S$ for which $d_X(x,y)>S$ implies $u(x,y)=0$.
\end{thm}

\begin{prop}
\label{properaction}
Let $G$ be a countable group acting properly isometrically on a proper metric space $X$. Then for any basepoint $x_0$ in $X$, the map $\psi\colon g\mapsto gx_0$ is a uniform embedding, and any map $\phi\colon Gx_0 \to G$ such that $\psi\circ\phi=\id_{Gx_0}$ is a coarse inverse to $\psi$. If moreover $G$ acts cocompactly then $\psi$ is a coarse equivalence between $G$ and $X$.
\end{prop}

\begin{proof}
First we verify that $\psi$ is a coarse map. For $g,h$ in $G$, we have $d_X(gx_0,hx_0)=d_X(x_0,g^{-1}hx_0)$ as the action is isometric. Given $R$, if $d_G(g,h)\leq R$ then $g^{-1}h$ lies in $B_R(e)$ which, by properness of the metric $d_G$ is finite, so $d_X(x_0,g^{-1}hx_0)$ is bounded by some number $S$. Thus for all $R$ there exists $S$ such that $d_G(g,h)\leq R$ implies $d_X(\psi(g),\psi(h))\leq S$. Properness of $\psi$ follows from properness of the action, so $\psi$ is a coarse map.

Now let $\phi\colon Gx_0 \to G$ be a splitting of $\psi$. Then for $x,y$ in $Gx_0$ we have $x=\phi(x)x_0,y=\phi(y)x_0$. If $d_X(x,y)\leq R$ then $d_X(x_0,\phi(x)^{-1}\phi(y)x_0)\leq R$. Properness of the action ensures that there are only finitely many elements $g$ of $G$ with $d_X(x_0,gx_0)\leq R$, so $\phi(x)^{-1}\phi(y)$ lies in $B_S(e)$, for some $S$, i.e.\ $d_G(\phi(x),\phi(y))\leq S$. The map $\phi$ is injective, thus it is proper, so $\phi$ is also a coarse map.

The composition $\psi\circ\phi$ is the identity on $Gx_0$ by definition, while $g^{-1}(\phi\circ\psi(g))$ is an element of $G$ fixing $x_0$. The stabiliser of $x_0$ is finite, so $d_G(g,\phi\circ\psi(g))$ is bounded as a function of $g$, i.e.\ $\phi\circ\psi$ is close to the identity on $G$. Thus $\psi$ is a uniform embedding from $G$ to $X$, and $\psi,\phi$ implement a coarse equivalence between $G$ and its image.

If the action is cocompact then for some $R$, the $G$ translations of the $R$-ball about $x_0$ cover $X$, so $Gx_0$ is $R$-dense in $X$. It follows immediately that the inclusion of $Gx_0$ into $X$ is a coarse equivalence, so $\psi \colon G \to X$ is a coarse equivalence.
\end{proof}

\begin{defn}
A kernel $k\colon X\times X \to \C$ has \emph{finite propagation} if there exists $R\geq 0$ such that $k(x,y)=0$ for $d(x,y)>R$. The $\emph{propagation}$ of $k$ is the smallest such $R$.
\end{defn}

If $X$ is a proper discrete metric space, and $k$ is a finite propagation kernel then for each $x$ there are only finitely many $y$ with $k(x,y)\neq 0$. Thus $k$ defines a linear map from $l^2(X)$ to itself, $k*v(x)=\sum_{y\in X} k(x,y)v(y)$. These linear maps are also said to have finite propagation. Note that if additionally $X$ has bounded geometry, then every bounded finite propagation kernel gives rise to a \emph{bounded} operator on $l^2(X)$.

\begin{defn}
The \emph{uniform Roe algebra}, $C^*_u(X)$, is the $C^*$-algebra completion of the algebra of bounded operators on $l^2(X)$ having finite propagation.
\end{defn}

\begin{defn}
A \emph{partial translation} of $X$ is a subset $t$ of $X\times X$ such that the coordinate projections of $t$ onto $X$ are injective, and such that $d_X(x,y)$ is bounded for $(x,y)\in t$.
\end{defn}
A partial translation can be viewed as a partially defined map from $X$ to $X$ which is close to the identity (where defined). Therefore a partial translation gives rise to a partial isometry of $l^2(X)$ which has finite propagation and hence is an element of $C^*_u(X)$. By definition, the partial isometry is defined by $t(\delta_y)=\delta_x$ if $t(y)=x$, and $t(\delta_y)=0$ if $t(y)$ is undefined.






\section*{Proof of Theorem \ref{mainthm}}
Fix a basepoint $x_0$ in $X$, and let $Y$ be the orbit of $x_0$ under the action of $G$. Given a finite propagation operator on $l^2(Y)$ we can extend it to $l^2(X)=l^2(Y)\oplus l^2(X\setminus Y)$ by defining it to be zero on $l^2(X\setminus Y)$. Thus the algebra $C^*_u(Y)$ is a subalgebra of $C^*_u(X)$, hence exactness for $C^*_u(X)$ implies exactness for $C^*_u(Y)$.

For each $y\in Y$, pick an element $\phi(y)$ of $G$ such that $\phi(y)x_0=y$. We construct a partially defined action of the group $G$ on the space $Y$ as follows. For $g\in G,y\in Y$ we define $g\act y=x$ if and only if $\phi(y)g^{-1}=\phi(x)$. The element $g\act y$ in $Y$ is uniquely determined if it exists, since then $g\act y=\phi(y)g^{-1}x_0$. However it will be undefined if $\phi(y)g^{-1}$ is not in the image of $\phi$. Note that $g\act $, viewed as a partially defined map from $X$ to $X$ is a partial translation, since
$$d_G(\phi(y),\phi(g\act y))=d_G(\phi(y),\phi(y)g^{-1})=d_G(e,g^{-1})$$
which is independent of $y$, and $\phi$ is a coarse equivalence by Proposition \ref{properaction}. Note that in the case where $G$ acts \emph{freely}, the map $\phi$ is uniquely determined and is a bijection between $Y$ and $G$. Using this to identify $Y$ with $G$, the action $g\act$ is identified with the action of $G$ on itself by right multiplication by $g^{-1}$.

Given $R,\eps$ we will produce a positive kernel $u$ on $Y$ with $(R,\eps)$-variation, and satisfying the hypothesis that $u(x,y)$ vanishes for $x,y$ far apart. Let $E_R$ be the set of elements of $G$ of the form $\phi(x)^{-1}\phi(y)$, $x,y$ in $Y$ with $d_X(x,y)\leq R$. As $\phi$ is a coarse map, there exists $S$ such that if $d_X(x,y)\leq R$ then $d_G(\phi(x),\phi(y))\leq S$, so $E_R$ is contained in the ball $B_S(e)$ in $G$, in particular it is finite. Elements of $G$ act as partial translations on $Y$ via $\act$, hence we can identify $E_R$ with a finite subset of $C^*_u(Y)$.

Using the characterisation of exactness (\cite{Oz} Lemma 2), there exists a completely positive finite rank map $\theta \colon C^*_u(Y) \to B(l^2(Y))$ such that:
\begin{enumerate}
\item $\theta$ has the form $\theta(\centerdot)=\sum_{i=1}^d \langle \delta_{a_i},\centerdot\delta_{b_i}\rangle T_i$, for some $a_i,b_i$ in $Y$, and $T_i$ in $B(l^2(Y))$, and
\item for all $g$ in $E_R$ we have $\|\theta(g)-g\|<\eps.$
\label{convergence}
\end{enumerate}
Let $F_R$ denote the set $\{a_i,b_i : i=1,\dots,d\}$, and note that for any partial translation $t$ such that the image of $F_R$ under $t$ does not meet $F_R$, we have $\theta(t)=0$. We now define
$$u(x,y)=\langle\delta_x,\theta(\phi(x)^{-1}\phi(y))\delta_y\rangle\text{ for $x,y$ in $Y$},$$
where the elements $\phi(x)^{-1}\phi(y)$ of $G$ are regarded as elements of $C^*_u(Y)$ as above.

First we will verify positivity of $u$; this is not immediate because as an operator on $l^2(X)$, $\phi(x)^{-1}\phi(y)$ is not necessarily the composition of the operators corresponding to $\phi(x)^{-1}$ and $\phi(y)$. This is because $(\phi(x)^{-1}\phi(y))\act y'$ may be defined even when $\phi(x)^{-1}\act(\phi(y)\act y')$ is not.

For each $y$ in $Y$, we define a bounded linear map $s_y$ from $l^2(Y)$ to $l^2(G)$ as follows. Let $s_y(\delta_{y'})=\delta_g$ where $g=\phi(y')\phi(y)^{-1}$ in $G$. Then its adjoint $s_y^*$ satisfies $s_y^*(\delta_g)=\delta_{x'}$ if there exists $x'$ with $\phi(x')\phi(x)^{-1}=g$, (such an $x'$ must be unique by injectivity of $\phi$) and is zero otherwise. Thus for $x,y,y'$ in $Y$, the vector $s_x^*s_y(\delta_{y'})$ is $\delta_{x'}$ for some $x'\in Y$ if we have
$$\phi(y')\phi(y)^{-1}=\phi(x')\phi(x)^{-1},$$
and it is zero otherwise. Note that we can rewrite this as $\phi(y')(\phi(x)^{-1}\phi(y))^{-1}=\phi(x')$
i.e.\ $x'=(\phi(x)^{-1}\phi(y))\act y'$. We conclude that
\begin{align*}
s_x^*s_y(\delta_{y'})&=\delta_{x'}=(\phi(x)^{-1}\phi(y))(\delta_{y'}) &&\text{ if $x'=(\phi(x)^{-1}\phi(y))\act y'$,}\\
s_x^*s_y(\delta_{y'})&=0=(\phi(x)^{-1}\phi(y))(\delta_{y'}) &&\text{ if $(\phi(x)^{-1}\phi(y))\act y'$ is undefined.}
\end{align*}
Hence $s_x^*s_y=\phi(x)^{-1}\phi(y)$ as an operator on $l^2(Y)$. The operators $(\phi(x)^{-1}\phi(y))=(s_x^*s_y)$ therefore form a positive matrix over $Y$, so positivity of $u$ follows from complete positivity of $\theta$.

We will now show that $u$ has $(R,\eps)$-variation. For $x,y$ with $d_X(x,y)\leq R$ we have $\phi(x)^{-1}\phi(y)$ in $E_R$, hence it follows from (\ref{convergence}) above that
$$\|\phi(x)^{-1}\phi(y)-\theta(\phi(x)^{-1}\phi(y))\|<\eps .$$
We recall that $(\phi(x)^{-1}\phi(y))\act y'$ is defined if $\phi(y')(\phi(x)^{-1}\phi(y))^{-1}$ is in the image of $\phi$. In particular for $y'=y$ we have that $\phi(y)(\phi(x)^{-1}\phi(y))^{-1}=\phi(x)$ which is in the image of $\phi$, so  $(\phi(x)^{-1}\phi(y))\act y$ is defined and equals $x$. Thus $\phi(x)^{-1}\phi(y)(\delta_y)=\delta_x$, so $\langle\delta_x,\phi(x)^{-1}\phi(y)(\delta_y)\rangle=1$. Hence
$$|1-u(x,y)|=\left|\left\langle\delta_x,\left(\phi(x)^{-1}\phi(y)-\theta\left(\phi(x)^{-1}\phi(y)\right)\right)(\delta_y)\right\rangle\right|<\eps.$$
We conclude that $u$ has $(R,\eps)$-variation as required.

Finally we will show that $u(x,y)$ vanishes for $d_X(x,y)$ sufficiently large. As $\phi$ is a uniform embedding, for all $R'$ there exists $S'$ such that if $d_G(\phi(x),\phi(y))\leq R'$ then $d_X(x,y)\leq S'$. As $F_R$ is a finite subset of $Y$ we note that $\{\phi(x')^{-1}\phi(y'):x',y'\in F_R\}$ is a finite subset of $G$, and choose $R'$ such that this is contained in the ball $B_{R'}(e)$ in $G$. Now if $x,y\in Y$ and $x',y'\in F_R$ with $(\phi(x)^{-1}\phi(y))\act y'=x'$ then $\phi(y')(\phi(x)^{-1}\phi(y))^{-1}=\phi(x')$ so $\phi(x')^{-1}\phi(y')=\phi(x)^{-1}\phi(y)$. Hence
$$d_G(e,\phi(x)^{-1}\phi(y))=d_G(e,\phi(x')^{-1}\phi(y'))\leq R',$$
i.e. $d_G(\phi(x),\phi(y))\leq R'$, so $d_X(x,y)\leq S'$. Thus if $d_X(x,y)>S'$ then for $y'$ in $F_R$, $(\phi(x)^{-1}\phi(y))\act y'$ cannot be an element of $F_R$. It follows that if $d_X(x,y)>S'$ then $\langle\delta_{x'},(\phi(x)^{-1}\phi(y))\delta_{y'}\rangle=0$ for all $x',y'$ in $F_R$, so $\theta(\phi(x)^{-1}\phi(y))$ vanishes. Hence $u(x,y)$ also vanishes as required.

We have shown that for all $R,\eps$ there exists a positive kernel on $Y$ with $(R,\eps)$-variation, and which vanishes for $x,y$ far apart. It follows by Theorem \ref{kernels} that $Y$ has property A. As $G$ is coarsely equivalent to $Y$, it too has property A, so by Ozawa's result $G$ is exact.

\end{document}